\documentclass[a4paper,11pt]{article}

\usepackage[T1]{fontenc}

\usepackage{amsmath,amsthm,amsfonts}
\usepackage{graphicx}

\renewcommand{\mid}{:}
\newcommand{\srb}{\textsc{srb}}
\newcommand{\N}{\mathbb{N}}
\newcommand{\Z}{\mathbb{Z}}

\newcommand{\R}{\mathbb{R}}

\newcommand{\stable}{\mathrm{s}}
\newcommand{\unstable}{\mathrm{u}}
\newcommand{\diff}[1]{\mathrm{d} #1}
\newcommand{\norm}[1]{\lVert #1 \rVert}

\newcommand{\seq}[1]{\underline{#1}}
\newcommand{\dimH}{\mathrm{dim}_\mathrm{H}\, }
\newcommand{\diam}{\mathrm{diam}\, }
\newcommand{\supp}{\mathrm{supp}\, }

\newtheorem{theorem}{Theorem}
\newtheorem{definition}{Definition}
\newtheorem{corollary}{Corollary}
\newtheorem{remark}{Remark}
\newtheorem{lemma}{Lemma}

\begin{document}
\author{Tomas Persson\footnote{Institute of Mathematics, Polish Academy of Sciences, ulica \'Sniadeckich~8, 00-950 Warszawa, Poland, \texttt{tomasp@impan.gov.pl}}}
\date{\today}
\title{\bf On the Hausdorff Dimension of\\ Piecewise Hyperbolic Attractors}

\maketitle

\begin{abstract}
We study non-invertible piecewise hyperbolic maps in the plane. The Hausdorff dimension of the attractor is calculated in terms of the Lyapunov exponents, provided that the map satisfies a transversality condition. Explicit examples of maps for which this condition holds are given.
\end{abstract}

\maketitle

\section{Introduction}

A general class of piecewise hyperbolic maps was studied by Pesin in \cite{Pesin}. Pesin proved the existence of \srb-measures and investigated their ergodic properties. Results from Pesin's article and Sataev's article \cite{Sataev} are described in Section \ref{sec:theory}. The assumptions in \cite{Pesin} and \cite{Sataev} did not allow overlaps of the images. Schmeling and Troubetzkoy extended in \cite{Schmeling-Troubetzkoy} the theory in \cite{Pesin} to allow maps with overlaping images.

Using the results of Pesin and techniques from Solomyak's paper \cite{Sol}, the author of this paper proved in \cite{Jag1} and \cite{Jag2} that for two classes of piecewise affine hyperbolic maps, there exists, for almost all parameters, an invariant measure that is absolutely continuous with respect to Lebesgue measure, provided that the map expands area. This result had previously been obtained for fat baker's transformations by Alexander and Yorke in \cite{Ale-Yor}.
The main difficulty that arises for the class of maps in \cite{Jag2} is that in difference from the fat baker's transformation the symbolic space associated to the systems, changes with the parameters, and also the \srb-measure changes in a way that is hard to control. By embedding all symbolic spaces into a larger space it was possible get sufficient control to prove the result.

Solomyak's proof in \cite{Sol} uses a transversality property of power series. The proofs in \cite{Jag1} and \cite{Jag2} uses that iterates of points can be written as power series with such a transversality property. For the possibility of writing iterates as power series, it is important that the directions of contraction is maped onto each other throughout the manifold. The method in \cite{Jag1} and \cite{Jag2} is therefore not good for proving similar results for more general maps. It should also be noted that this method only gives results that holds for almost every map, with respect to some parameter.

Tsujii studied in \cite{Tsujii_fat} a class of area-expanding solenoidal attractors and proved that generically these systems has an invariant measure that is absolutely continuous with respect to Lebesgue measure. Tsujii also used a transversality condition, but in a different way. Instead of transversality of power series, Tsujii used transversality of intersections of iterates of curves. This technique makes it possible to show the existence of an absolutely continuous invariant measure for a fixed system, provided that the appropriate transversality condition is satisfied. Tsujii proved that this transversality condition is generically satisfied.

In this paper we will use this idea from Tsujii's article \cite{Tsujii_fat} to prove a formula for the dimension of the attractor for some piecewise hyperbolic maps in the plane, provided that a transversality condition is satisfied. This is done by an estimate of the dimension from below. This estimate coinsides with a previously known estimate from above (see \cite{Falconer2} and \cite{Schmeling-Troubetzkoy}) and thereby provide the following formula for the dimension:
\[
  \dimH \Lambda = 1 - \frac{\chi_\unstable}{\chi_\stable},
\]
where $\Lambda$ denotes the attractor and $\chi_\unstable$ and $\chi_\stable$ denote the positive and the negative Lyapunov exponents.
This formula has previously been proved by Falconer in \cite{Falconer2} and by Simon in \cite{Simon}, but for a much smaller class of systems. Both Falconer and Simon considered maps that are scew-products with the underlying shift being a full shift on $n$ symbols. These restrictions are not assumed in this paper.
Hence, this paper generalises the results of Falconer and Simon.

We will also need the assumption that the multiplicity entropy is zero, which also is the case in Falconer's and Simon's results. This seems often to be the case, and we provide a condition which guaranties that the multiplicity entropy is zero.

\section{Outline of the Paper}

In Section~\ref{sec:theory} we present the general theory of piecewise hyperbolic maps, that will be used later in the paper. In Section~\ref{sec:trans} we introduce a transversality condition. Under the assumption that this transversality condition holds, a theorem that estimates the dimension from below is stated in Section~\ref{sec:dimension}. This estimate gives the dimension formula. The theorem is proved in Section~\ref{sec:proof} and Section~\ref{sec:example} contains explicit examples of maps that satisfy the assumptions of this theorem. There are also examples that the dimension formula may fail if the transversality condition does not hold.

\section{Piecewise Hyperbolic Maps} \label{sec:theory}

There is a study of general piecewise hyperbolic maps in Pesin's article \cite{Pesin}. He studied maps of the following form.

Let $M$ be a smooth Riemannian manifold with metic $d$,
let $K \subset M$ be an open, 
bounded and connected set and let $N \subset K$ be a closed set in $K$. The set $N$ is called the discontinuity set. Let $f \colon K \setminus N \rightarrow K$.
 
Put
\begin{align*}
  &K^+ = \{\,x \in K \mid f^n(x) \not \in N \cup \partial K, \
  n=0,1,2,\ldots\,\}, \\
  &D = \bigcap_{n \in \N} f^n (K^+).
\end{align*}
The attractor of $f$ is the set $\Lambda = \overline{D}$.

The maps studied in \cite{Pesin} were assumed to satisfy the
following conditions.
\begin{equation*} \tag{A1} \label{A1}
\begin{minipage}[t]{0.8\textwidth} 
$f \colon K \setminus N \rightarrow f(K \setminus N)$ is a
    $C^2$-diffeomorphism.
\end{minipage}
\end{equation*}

Let $N^+ = N \cup \partial K$ and
\begin{multline*}
       N^- = \{\,y \in K \mid \exists z_n \in K \setminus N^+,\ z \in N^+\\ \mathrm{such}\ \mathrm{that}\ y = f(z),\ z_n \rightarrow z,\ f(z_n) \rightarrow y\,\}.
\end{multline*}
One might want to think of $N^-$ as the image of $N^+$ although $f$ is not
defined on $N^+$.
We can now formulate the second and third assumption.

\begin{equation*} \tag{A2} \label{A2}
  \begin{minipage}[t]{0.8\textwidth} 
    There exists $C > 0$ and $\alpha \geq 0$ such that
     \begin{align*}
       \norm{\diff{}^2_x f} &\leq C d(x, N^+)^{-\alpha},
       &&\forall x \in K\setminus N,\\
       \norm{\diff{}^2_x (f^{-1})} &\leq C d(x, N^-)^{-\alpha},      &&\forall x \in f(K\setminus N).
     \end{align*}
  \end{minipage}
\end{equation*}

\begin{equation*} \tag{A3} \label{A3}
\begin{minipage}[t]{0.8\textwidth} 
For $\varepsilon > 0$ and $l=1,2,\ldots$ let
    \begin{align*}
      D_{\varepsilon,l}^+ &=\{\,x \in K^+ \mid d(f^n (x), N^+) \geq
      l^{-1} e^{-\varepsilon n},\ n\in\N\,\},\\
      D_{\varepsilon,l}^- &=\{\,x \in \Lambda \mid d(f^{-n} (x), N^-) \geq
      l^{-1} e^{-\varepsilon n},\ n\in\N\,\},\\
      D_{\varepsilon, l}^0 &= D_{\varepsilon, l}^+ \cap D_{\varepsilon, l}^-, \\
      D_\varepsilon^0 &= \bigcup_{l \geq 1} 
        (D_{\varepsilon, l}^+ \cap D_{\varepsilon, l}^-).
    \end{align*}
    The set $D_\varepsilon^0$ is not empty for sufficiently small
    $\varepsilon > 0$. (Here sufficiently small means so small that there are local unstable manifolds.)
\end{minipage}
\end{equation*}
The attractor is called regular if \eqref{A3} is satisfied. For a given map, it is usually not apperent whether the condition \eqref{A3} is satisfied or not. There exist however conditions that implies \eqref{A3} and are such that it easily can be checked if they hold true. These conditions are given in the end of this section.

\begin{equation*} \tag{A4} \label{A4}
\begin{minipage}[t]{0.8\textwidth} 
There exists $C > 0$ and $0 < \lambda < 1$ such that 
    for every $x \in K\setminus N^+$ there exists cones
    $C^\stable (x), C^\unstable (x) \subset T_x M$ such that the angle
    between $C^\stable (x)$ and $C^\unstable (x)$ is uniformly bounded
    away from zero,
    \begin{align*}
      &\diff{}_x f(C^\unstable (x)) \subset C^\unstable (f(x))
      && \forall x \in K\setminus N^+, \\
      &\diff{}_x (f^{-1}) (C^\stable (x)) \subset C^\stable (f^{-1}(x))
      && \forall x \in f(K\setminus N^+),
    \end{align*}
    and for any $n > 0$
    \begin{align*}
      &\norm{\diff{}_x f^n (v)} \geq C \lambda^{-n} \norm{v},
      && \forall x \in K^+,\ \forall v \in C^\unstable (x),\\
      &\norm{\diff{}_x f^{-n} (v)} \geq C \lambda^{-n} \norm{v},
      && \forall x \in f^n(K^+),\ \forall v \in C^\stable (x).
    \end{align*}
\end{minipage}
\end{equation*}

The last assumption makes it possible to define stable and unstable manifolds, $W^\stable (x)$ and $W^\unstable (x)$ as well as local ones
for any $x \in D_\varepsilon^0$.

The condition
\begin{equation*} \tag{A3$'$} \label{A3b}
\begin{minipage}[t]{0.8\textwidth} 
There exists a point $x \in D_\varepsilon^0$ and
    $C,t,\delta_0 > 0$
    such that for any $0 < \delta < \delta_0$ and any
    $n \geq 0$
    \[
      \nu^\unstable (f^{-n} (U(\delta, N^+))) < C \delta^t,
    \]
    where $\nu^\unstable $ is the measure on the local unstable manifold
    of $x$, induced by the Riemannian measure, and $U(\delta, N^+)$ is an open $\delta$-neigbourhood of $N^+$.
\end{minipage}
\end{equation*}
implies condition \eqref{A3}.
Pesin proved the following theorem.

\begin{theorem}[{\rm Pesin \cite{Pesin}}] \label{the:pesin}
  Assume that $f$ satisfies the assumptions \eqref{A1}--\eqref{A4} and \eqref{A3b}.
  Then there exists an $f$-invariant measure $\mu$ such that
  $\Lambda$ can be decomposed $\Lambda = \bigcup_{i\in\N} \Lambda_i$
  where
  \begin{itemize}
    \item $\Lambda_i \cap \Lambda_j = \emptyset$, if $i \ne j$,
    \item $\mu (\Lambda_0) = 0$, $\mu (\Lambda_i) > 0$  if  $i > 0$,
    \item $f(\Lambda_i) = \Lambda_i,\ f|_{\Lambda_i}$ is ergodic,
    \item for $i > 0$ there exists $n_i > 0$ such that 
          $(f^{n_i} |_{\Lambda_i}, \mu)$ is isomorphic to a Bernoulli shift.
  \end{itemize}
  The metric entropy satisfy
  \[
    h_\mu (f) = \int \sum \chi_i (x) \: \diff{\mu} (x),
  \]
  where the sum is over the positive Lyapunov exponents $\chi_i (x)$.
\end{theorem}

The measure $\mu$ in Theorem \ref{the:pesin} is called \srb-measure
(or Gibbs u-measure). For piecewise hyperbolic maps the \srb-measures
are characterised by the property that their conditional measures on unstable manifolds are absolutely continuous with respect to Lebesgue measure and the set of typical points has positive Lebesgue measure.

For a somewhat smaller class of maps Sataev proved in \cite{Sataev} that
the ergodic com\-ponents of the \srb-measure (the sets $\Lambda_i$ in
Theorem \ref{the:pesin}) are finitely many.

\subsection{Non-Invertible Piecewise Hyperbolic Maps} \label{ssec:noninv}

The maps studied by Pesin and Sataev are all invertible on their images.
Schmeling and Troubetzkoy generalised in \cite{Schmeling-Troubetzkoy} the results of Pesin to non-invert\-ible maps: If 
\begin{equation*} \tag{A5} \label{A5}
\begin{minipage}[t]{0.8\textwidth} 
the set
$K \setminus N$ can be decomposed into finitely many sets $K_i$ such that
$f \colon K_i \rightarrow f(K_i)$ can be extended to a diffeomorphism from $\overline{K}_i$ to $\overline{f(K_i)}$
\end{minipage}
\end{equation*}
and $f$ satisfies the
assumptions \eqref{A2}--\eqref{A4} and \eqref{A3b}, then the statement of Theorem~\ref{the:pesin} is still valid.
Note that $f(K_i) \cap f(K_j)$ is allowed to be non-empty so that
$f \colon K\setminus N \rightarrow f(K\setminus N)$ is not a diffeomorphism. Schmeling and Troubetzkoy proved their result by lifting the
map and the set $K$ to a higher dimension; Let $\hat{K} = K \times [0,1]$,
$\hat{K}_i = K_i \times [0,1]$ and
\[
  \hat{f} |_{K_i} \colon (x,t) \mapsto (f(x), \tau t + i/p), \quad i=0,1,\ldots, p-1,
\]
where $\tau < 1$ and $p$ is the number of sets $K_i$. The map $\hat{f}$ is then invertible if $\tau$ is sufficiently small and then $\hat{f}$ satisfies the assumptions of Theorem \ref{the:pesin}, in particular there is an \srb-measure $\hat{\mu}$ on the lifted set $\hat{K}$. The projection of this measure to the set $K$ was shown to be an \srb-measure of the original map $f$, in the sence that the set of typical points with respect to the projected measure has positive Lebesgue measure.

We will let $\hat{D}$, $\hat{D}_{\varepsilon,l}^0, \ldots$ denote the lifted variants of the corresponding sets $D$, $D_{\varepsilon,l}^0, \ldots$

It is often hard to check whether \eqref{A3b} holds. It is proved in \cite{Schmeling-Troubetzkoy} that if $f$ satisfies \eqref{A2}, \eqref{A4}, \eqref{A5} and 
the assumptions \eqref{A6}--\eqref{A8} below, then $f$ satisfies condition \eqref{A3b}, and hence also \eqref{A3}.
\begin{equation*} \tag{A6} \label{A6}
\begin{minipage}[t]{0.8\textwidth} 
  The sets $\partial K$ and $N$ are unions of finitely many smooth curves such that the angle between these curves and the unstable cones are bounded away from zero.
\end{minipage}
\end{equation*}
\smallskip
\begin{equation*} \tag{A7} \label{A7}
\begin{minipage}[t]{0.8\textwidth} 
  The cone families $C^\unstable (x)$ and $C^\stable (x)$ depends continuously on $x \in K_i$ and they can be extend continuously to the boundary.
\end{minipage}
\end{equation*}
\smallskip
\begin{equation*} \tag{A8} \label{A8}
\begin{minipage}[t]{0.8\textwidth} 
  There is a natural number $q$ such that at most $L$ singularity curves of $f^q$ meet at any point, and $a^q > L + 1$ where
\[
  a = \inf_{x \in K \setminus N} \inf_{v \in C^\unstable (x)} \frac{|\diff{}_x f (v)|}{|v|}.
\]

\end{minipage}
\end{equation*}

\subsection{Multiplicity entropy} \label{ssec:multiplicitydefinition}

The assumption \eqref{A8} implies that the multiplicity entropy \cite{buzzi} is not larger than $\log (L+1)$. We will need a stronger assumption (assumption \eqref{A9} below) than \eqref{A8}, namely that the multiplicity entropy is zero. We will show that this is satisfied under rather mild assumptions on the map. Let us start with defining the multiplicity entropy, and then give the assumption \eqref{A9}.

Let $\mathcal{K} = \{K_1, \ldots, K_p\}$ be the partion of $K$ into sets on which $f$ is continuous, and let $\mathcal{K}_n$ be the corresponding partition for the map $f^n$. Let $k_n$ be the maximal numbers of elements of $\mathcal{K}_n$ that meet in one point. The multiplcity entropy $h_\mathrm{mult} (f)$ is defined as
\[
  h_\mathrm{mult} (f) = \limsup_{n \to \infty} \frac{1}{n} \log k_n.
\]

We can now give our next assumption.
\begin{equation*} \tag{A9} \label{A9}
  h_\mathrm{mult} (f) = 0.
\end{equation*}
One might wounder how general this condition is. The author of this paper knows of no example of a map in dimension two, satisfying \eqref{A1}--\eqref{A8}, with one positive and one negative Lyapunov exponent and such that \eqref{A9} is not satisfied. In Section~\ref{sec:multiplicity} we give sufficient conditions for the map to satisfy  \eqref{A9}, and hence also \eqref{A8}. These conditions are for instance satisfied by Belykh maps. Hence Belykh maps have zero multiplicity entropy.

For future use, we note that condition \eqref{A9} implies that the topological entropy is equal to the entropy of the \srb-measure. This follows by the result of Kruglikov and Rypdal in \cite{KruglikovRypdal}, that $h_\mathrm{top} \leq \chi_\unstable + h_\mathrm{mult}$ (in the case of a map on the plane with one positiv and one negative Lyaponov exponent; the statement in \cite{KruglikovRypdal} is for any dimension).

\section{A Transversality Condition} \label{sec:trans}

Let $\varepsilon > 0$ and $0 < \delta < 1$. We will say that an intersection of two smooth curves $\gamma_1$ and $\gamma_2$ is $(\varepsilon, \delta)$-transversal if for any ball $B_\varepsilon$ of radius $\varepsilon$ intersecting both $\gamma_1$ and $\gamma_2$, there exist points $x_1 \in B_\varepsilon \cap \gamma_1$ and $x_2 \in B_\varepsilon \cap \gamma_2$ such that the following holds true. If $d_1$ and $d_2$ are the induced metrics on $\gamma_1$ and $\gamma_2$ respectively, then the intersection of the open sets
\begin{equation} \label{transversal}
  \bigcup_{y \in \gamma_i \cap B_\varepsilon} B (y, \delta d_i (x_i, y)), \quad i = 1,2,
\end{equation}
is empty. The symbols $B (x, r)$ denotes the open ball of radius $r$ around $x$. 
Note that if $\gamma_1$ and $\gamma_2$ intersect $(\varepsilon, \delta)$-transversal then the intersection $\gamma_1 \cap \gamma_2$ can be empty.

\begin{definition} \label{def:trans}
We will say that a piecewise hyperbolic system $f \colon K \setminus N \to K$ satisfies condition \eqref{T} if
\begin{equation*} \tag{T} \label{T}
\begin{minipage}[t]{0.8\textwidth} 
there exists numbers $\varepsilon, \delta > 0$ such that if\/ $\gamma_1$ and $\gamma_2$ are two smooth curves such that every tangent lies in the unstable cone, and $\gamma_1$ and $\gamma_2$ are in different $K_i$, then the curves $f (\gamma_1)$ and $f (\gamma_2)$ intersect $(\varepsilon, \delta)$-transversal.
\end{minipage}
\end{equation*}

\end{definition}

\section{Dimension of the Attractor} \label{sec:dimension}

Consider a map $f \colon K \setminus N \to K \subset \R^2$ that satisfies the conditions \eqref{A2}, \eqref{A4} and \eqref{A5}--\eqref{A9}. We denote by $\chi_\stable (x) < 0 < \chi_\unstable (x)$ the two Lyapunov exponents at the point $x$ if they exist. If $\Lambda_1$ is an ergodic component of the attractor, then the Lyapunov exponents are constant almost everywhere with respect to the \srb-measure on $\Lambda_1$, and we write $\chi_\stable (x) = \chi_\stable$ and $\chi_\unstable (x) = \chi_\unstable$ for almost every $x$.

\begin{theorem} \label{the:dimension}
  Suppose that $f \colon K\setminus N \to K \subset \R^2$ is a piecewise hyperbolic map that satisfies the conditions \eqref{T}, \eqref{A2}, \eqref{A4} and \eqref{A5}--\eqref{A9}. Let $\Lambda_1$ be an ergodic component of\/ the attractor, whith one positive and one negative Lyapunov exponent. Then the Hausdorff dimension of\/ $\Lambda_1$ satisfies
  \[
     \dimH \Lambda_1 \geq \min \biggl\{2,\, 1 - \frac{\chi_\unstable}{\chi_\stable} \biggr\}.
  \]
\end{theorem}

Theorem \ref{the:dimension} is proved in Section \ref{sec:proof}.

Note that in \cite{Schmeling-Troubetzkoy}, it is proved that $\dimH \Lambda_1 \leq 1 - \chi_\mathrm{u} / \chi_\mathrm{s}$ with equality if and only if $f$ restricted to $\Lambda_1$ is almost everywhere invertible, meaning that $f$ is invertible on a set of full measure. Hence we get the following corollary.

\begin{corollary} \label{cor:eq}
  If the assumptions of\/ Theorem~\ref{the:dimension} are satisfied then
  \[
    \dimH \Lambda_1 = \min \biggl\{2,\, 1 - \frac{\chi_\mathrm{u}}{\chi_\mathrm{s}} \biggr\},
  \]
  and $f$ is invertible almost everywhere on $\Lambda_1$ if and only if\/ $\chi_\unstable + \chi_\stable \leq 0$.
\end{corollary}

\begin{remark}
In case the transverality condition \eqref{T} is not satisfied we can only give the trivial estimate $\dimH \Lambda_1 \geq 1$. Indeed, the map $f \colon [0,1]^2 \to [0,1]^2$ defined by $f \colon (x_1,x_2) \mapsto (x_1 /2, 2 x_2\ \mathrm{mod}\, 1)$ has the attractor $\Lambda_1 = \{\, (x_1,x_2) : x_1 = 0,\ 0\leq x_1 < 1 \,\}$, and so $\dimH \Lambda_1 = 1$. Moreover, any map satisfying the conditions \eqref{A2}, \eqref{A4} and \eqref{A5}--\eqref{A8} has an attractor $\Lambda_1$ that contains curves of unstable manifolds. This implies that $\dimH \Lambda_1 \geq 1$. So, unless one impose an additional condition, such as \eqref{T}, one can not get a better estimate than $\dimH \Lambda_1 \geq 1$.
\end{remark}

\begin{remark}
It should be noted that if $f$ satisfies the conditions in Theorem~\ref{the:dimension}, then so does any sufficiently small smooth perturbation of $f$.
\end{remark}

\section{Vanishing multipicity entropy} \label{sec:multiplicity}

In this section we give a condition which guaranties that the multiplicity entropy is zero.

\begin{theorem} \label{the:multiplicity}
  Let $K, N \subset \R^2$ where $N$ is a union of smooth curves, and let $f \colon K \setminus N \to K$ satisfy conditions \eqref{A1}--\eqref{A7}.
  Assume that there is a family of cones $C^\mathrm{d} (p, \gamma) \subset T_p \R^2$ where $p$ is a point on a smooth curve $\gamma \subset N$, such that
  \begin{equation} \label{multcones}
    C^\mathrm{d} (p, \gamma) \cap C^\unstable (p) = \{0\} \quad \mathrm{and} \quad \diff{f} (C^\mathrm{d}) \subset C^\unstable.
  \end{equation}
  Then the multiplicity entropy of $f$ is zero.
\end{theorem}

\begin{remark}
  The condition $C^\mathrm{d} (p, \gamma) \cap C^\unstable (p) = \{0\}$ in Theorem~\ref{the:multiplicity} is nonsence since $C^\unstable (p)$ is not defined for $p \in N$. But $C^\unstable (p)$ depends continuously on $p \in K_i$ so the condition should be understand as $C^\unstable (p)$ replaced by its limit for each $K_i$ that meet $p$.
\end{remark}

\begin{figure}
  \begin{center}
    \setlength{\unitlength}{0.008\textwidth}
    {\begin{picture}(100,75)
    \put(0,0){{\includegraphics[width=0.8\textwidth]{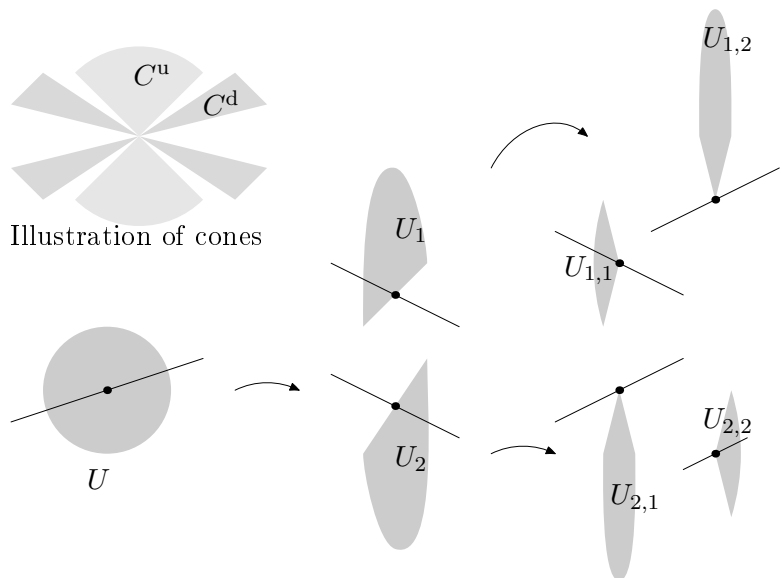}}}
    \put(10,12){$U$}
    \put(50,45){$U_1$}
    \put(50,15){$U_2$}
    \put(72,40){$U_{1,1}$}
    \put(90,70){$U_{1,2}$}
    \put(78,10){$U_{2,1}$}
    \put(90,20){$U_{2,2}$}
    \put(16,65){$C^\unstable$}
    \put(25,61){$C^\mathrm{d}$}
    \put(0,44){Illustration of cones}
    \end{picture}}
  \end{center}
  \caption{Illustration to the proof of Theorem~\ref{the:multiplicity}. Note that $U_{1,2}$ and $U_{2,1}$ cannot be cut throught $f^2 (p)$ since the slopes of the discontinuities are too small.} \label{fig:multiplicity}
\end{figure}

\begin{proof}
  For simplicity, let us start with the case that the curves of $N$ do not intersect. Let $p \in N$. We will iterate $p$ and see in how many pieces a small neigborhood $U$ of $p$ is cut by a curve in $N$ that goes through $f^n (p)$. Of cource, $f^n (p)$ is not defined but we will use this notation for simplicity, for the collection of accumulation points of $f^n (q)$, when $q \to p$.

  In the first iterate $U$ is cut through $p$ in at most two pieces, which we denote by $U_{1}$ and $U_{2}$ (or just $U_1$ if $U$ is not cut). In the next iterate, each of the pieces $U_1$ and $U_2$ is cut throught $f(p)$ in at most two pieces. Denote by $U_{1,1}$ and $U_{1,2}$ the pieces of $U_1$ and similarly for $U_2$. 

  By the property \eqref{multcones}, one of $U_{1,1}$ and $U_{1,2}$ lies in the cone $C^\unstable (f^2(p))$ and no iterate of this piece will be cut through $f^n (p)$ for any $n$. The same argument holds for the pieces $U_{2,1}$ and $U_{2,2}$. So we now have at most four pieces of which at most two can be cut in future iterations. There is a picture of this in Figure~\ref{fig:multiplicity}.

  By induction we get that after $n$ iterates $f^n (U)$ consists of at most $2n$ pieces. This shows that the multiplicity entropy is zero.

  The case with $N$ containing curves that cut each other is similar. If at most $L$ curves meet in one point, we get that after $n$ iterates, $U$ consists of at most $2 (L + 1) n$ pieces of which at most two can be cut through $f^n (p)$.
\end{proof}

\section{An Example} \label{sec:example}

In this section we give an example of maps satifying the assumptions of Theorem~\ref{the:dimension}.

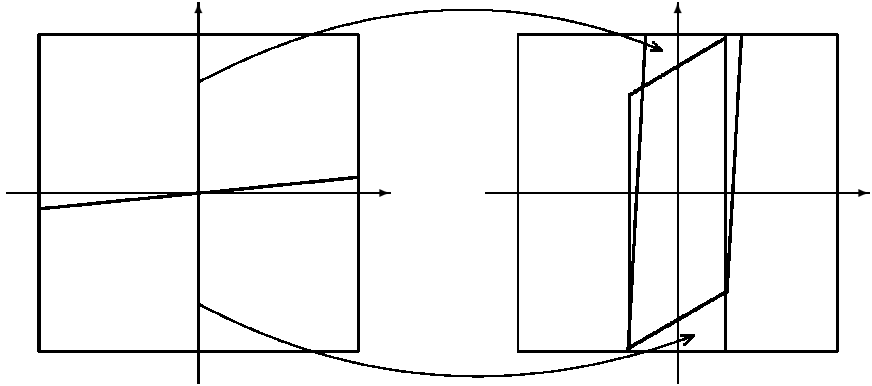
\begin{figure}
\begin{center}
\setlength{\unitlength}{0.166\textwidth}
\begin{picture}(6,3)(-1.5,-1.5)
\put(-1.2,0){\vector(1,0){2.4}}
\put(0,-1.2){\vector(0,1){2.4}}

\put(1.8,0){\vector(1,0){2.4}}
\put(3,-1.2){\vector(0,1){2.4}}

\thicklines
\put(-1,-1){\line(1,0){2}}
\put(-1,-1){\line(0,1){2}}
\put(1,1){\line(-1,0){2}}
\put(1,1){\line(0,-1){2}}
\qbezier(-1,-0.1)(0,0)(1,0.1)

\thinlines
\put(2,-1){\line(1,0){2}}
\put(2,-1){\line(0,1){2}}
\put(4,1){\line(-1,0){2}}
\put(4,1){\line(0,-1){2}}

\thicklines
\qbezier( 3.4, 1 )( 3.1, 1 )( 2.8, 1 )
\qbezier( 2.8, 1 )( 2.745, 0.01 )( 2.69, -0.98 )
\qbezier( 2.69, -0.98 )( 3, -0.8 )( 3.31, -0.62 )
\qbezier( 3.31, -0.62 )( 3.355, 0.19 )( 3.4, 1 )

\qbezier( 3.3, -1 )( 3.3, -0.01 )( 3.3, 0.98 )
\qbezier( 3.3, 0.98 )( 3, 0.8 )( 2.7, 0.62 )
\qbezier( 2.7, 0.62 )( 2.7, -0.19 )( 2.7, -1 )
\qbezier( 2.7, -1 )( 2.95, -1 )( 3.3, -1 )

\thinlines
\qbezier(0,0.7)(1.5,1.5)(2.9,0.9)
\qbezier(2.9,0.9)(2.875,0.925)(2.85,0.95)
\qbezier(2.9,0.9)(2.85,0.9)(2.83,0.9)
\qbezier(0,-0.7)(1.5,-1.5)(3.1,-0.9)
\qbezier(3.1,-0.9)(3.05,-0.9)(3.02,-0.9)
\qbezier(3.1,-0.9)(3.05,-0.95)(3.04,-0.96)

\end{picture}

\caption{A picture of $f$ with $\rho = 0.1$, $\psi_1=\psi_2=0$, $\gamma =1.8$, $\lambda=0.3$, $k=0.1$,  $a_1=a_2=0$ and $-b_1 = b_2 = 0.8$} \label{fig:map}
\end{center}
\end{figure}

Let $K=(-1,1) \times (-1,1)$ be a square. Take $-1 < k < 1$ and let $N = \{\, (x_1, x_2) \in K \mid x_2 = k x_1 \,\}$ be the singularity set. Take $\rho \ne 0$ and let $\psi_1$ and $\psi_2$ be two $C^2$ functions, such that $|\psi_1'|, |\psi_2'| < \rho_\psi < |\rho|/2$.
We take parameters $\frac{1}{2} < \lambda < 1$, $1 < \gamma < 2$, $a_1$, $a_2$, $b_1$ and $b_2$ such that the map $f$ defined by
\begin{equation} \label{sqewbelykh}
  f(x_1,x_2) = \left\{ \begin{array}{lll} (\lambda x_1 + a_1 + \rho x_2 + \psi_1 (x_2), & \gamma x_2 + b_1 ) & \mathrm{if }\ x_2 > kx_1 \\
  (\lambda x_1 + a_2+ \psi_2 (x_2), &\gamma x_2 + b_2 ) & \mathrm{if }\ x_2 > kx_1 \end{array} \right.
\end{equation}
maps $K \setminus N$ into $K$.  There is a picture of $f$ in Figure~\ref{fig:map}.

The case $\rho \ne 0$, $k = \psi_1 = \psi_2 = 0$ and $\gamma = 2$ is threated by Falconer in \cite{Falconer2}. He proved that for almost all parameters $\gamma$ and $\lambda$, the dimension satisfies $\dimH \Lambda = 1 - \log \gamma / \log \lambda$. The case $k=0$ and $\gamma = 2$ is covered by Simon's paper \cite{Simon}. He proved equality for all parameters. We prove that we have equality for all parameters in  also when $k$, $\psi_1$ and $\psi_2$ are not nessesarily zero. 
More precisely, we use Theorem \ref{the:dimension} to prove the following theorem.

\begin{theorem} \label{the:belykh}
If $a_1$, $a_2$, $- b_1 = b_2 = (\gamma -1)$ and
  \[
    (\gamma, \lambda, k, \rho) \in \{\,(\gamma, \lambda, k, \rho) : \gamma > 2 \lambda,\  \rho \ne 0 \, \}
  \]
  are numbers such that $f \colon K \setminus N \to K$, then $f \colon K\setminus N \to K$ defined by \eqref{sqewbelykh} has an attractor $\Lambda$ with dimension
  \begin{equation} \label{belykhdim}
     \dimH \Lambda = \min \biggl\{2,\, 1 - \frac{\log \gamma}{\log \lambda} \biggr\}.
  \end{equation}
\end{theorem}

Let $\psi_1=\psi_2=0$, $1< \gamma<2$, $0<\lambda < 1$, $a_1=a_2=0$ and $b_1 = -b_2 = 1-\gamma$. Then if $\rho = 0$, the attractor is $\Lambda = \{\, (x_1,x_2) : x_1 = 0,\ |x_2| \leq \gamma-1 \,\}$, and so $\dimH \Lambda = 1$. If $\rho \ne 0$ and $\gamma > 2 \lambda$ then the dimension $\dimH \Lambda$ is given by \eqref{belykhdim}. The dimension can be made arbitrarily close to $2$ by choosing $\lambda$ close to $1$. Then the dimension is bounded away from 1 for any $\rho \ne 0$ but the dimension is 1 for $\rho = 0$.

\begin{proof}[Proof of Theorem \ref{the:belykh}.] 

It is clear from Theorem~\ref{the:multiplicity} that $f$ has zero multiplicity entropy if $k \ne 0$. If $k=0$ then the multiplicity entropy is trivialy zero.

We claim that if $\gamma > 2 \lambda$ and $\rho \ne 0$ then $f$ satisfies condition \eqref{T}. Let us prove this claim. It is clear that the cone spanned by the vectors 
\[
 \Bigl(\frac{-\rho_\psi}{\gamma - \lambda}, 1 \Bigr) \quad \mathrm{and}\ \quad \Bigl(  \frac{\rho + \rho_\psi}{\gamma - \lambda}, 1 \Bigr)
\]
defines an unstable cone family at any point of $K \setminus N$. Denote this cone by $C^\unstable$.

If $\sigma_1 \subset K \cap \{ x_2 > kx_1\}$ and $\sigma_2 \subset K \cap \{ x_2 < kx_1\}$ are two curves such that if $v_1$ and $v_2$ are two tangent vectors of the curves, then $v_1, v_2 \in C^\unstable$. The vectors $v_1$ and $v_2$ are mapped by $\diff{}_x f$ to
\[
  u_1 = \left[ \begin{array}{ll} \lambda & \rho + \psi_1(x_2) \\ 0 & \gamma \end{array} \right] v_1 \quad \mathrm{and}\ \quad u_2 = \left[ \begin{array}{ll} \lambda & \psi_2(x_2) \\ 0 & \gamma \end{array} \right] v_2
\]
respectively. One checks that $u_1$ is contained in the cone spanned by
\[
  \Bigl( - \rho_\psi \frac{\lambda}{\gamma(\gamma-\lambda)} + \frac{\rho - \rho_\psi}{\gamma}, 1 \Bigr) \quad \mathrm{and}\ \quad \Bigl( (\rho + \rho_\psi)  \frac{\lambda}{\gamma(\gamma-\lambda)} + \frac{\rho + \rho_\psi}{\gamma}, 1 \Bigr)
\]
and $u_2$ is contained in the cone spanned by 
\[
  \Bigl( - \rho_\psi \frac{\lambda}{\gamma(\gamma-\lambda)} + \frac{ - \rho_\psi}{\gamma}, 1 \Bigr) \quad \mathrm{and}\ \quad \Bigl( ( \rho + \rho_\psi )  \frac{\lambda}{\gamma(\gamma-\lambda)} + \frac{\rho_\psi}{\gamma}, 1 \Bigr)
\]
The intersection of these two cones is trivial if
\[
   - \rho_\psi \frac{\lambda}{\gamma(\gamma-\lambda)} + \frac{\rho - \rho_\psi}{\gamma} > ( \rho + \rho_\psi ) \frac{\lambda}{\gamma(\gamma-\lambda)} + \frac{\rho_\psi}{\gamma}, 
\]
or equivalently, if $\gamma > 2 \lambda$. This proves the claim.

By Corollary~\ref{cor:eq} it now follows that
\begin{equation*} 
 \dimH \Lambda = 1 - \frac{\log \gamma}{\log \lambda},
\end{equation*}
unless $\log \gamma + \log \lambda > 0$, in which case $\dimH \Lambda = 2$.
\end{proof}

\begin{figure}
\setlength{\unitlength}{0.166\textwidth}
\begin{minipage}[t]{0.49\textwidth}
\begin{center}
 \includegraphics[width=0.95\textwidth]{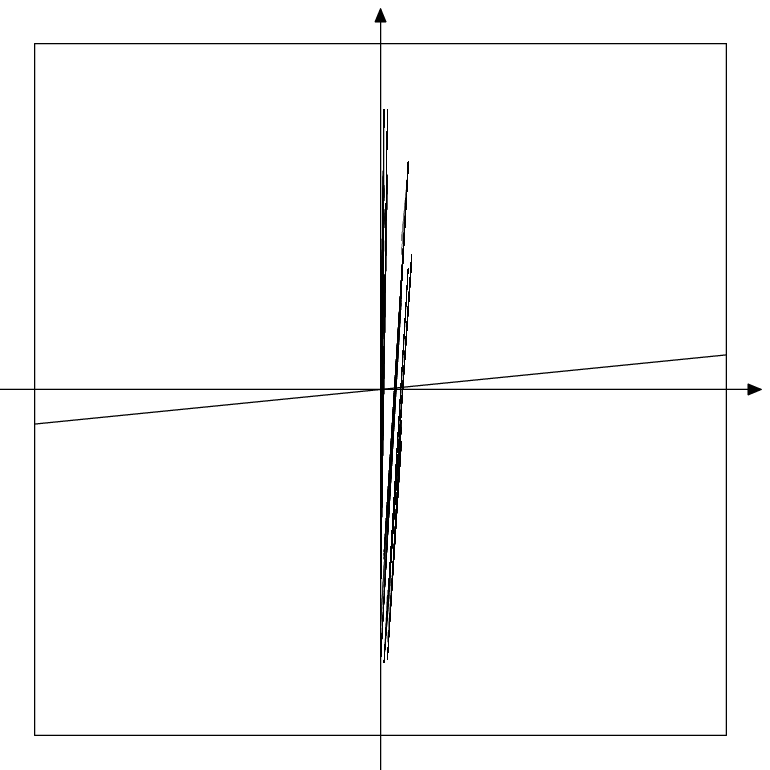}
\caption{The attractor $\Lambda$ of the map in Figure~\ref{fig:map}.} \label{fig:att}
\end{center}
\end{minipage}
\hspace{\stretch{1}}
\begin{minipage}[t]{0.49\textwidth}
\begin{center}
 \includegraphics[width=0.95\textwidth]{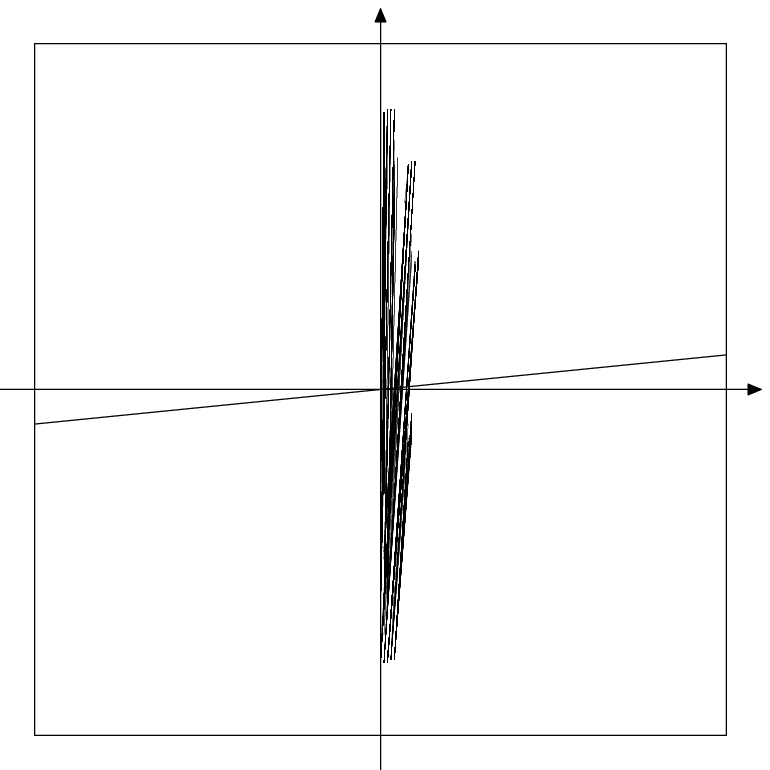}
\caption{The attractor $\Lambda$ of the map in Figure~\ref{fig:map}, but with $\lambda = 0.5$.} \label{fig:att2}
\end{center}
\end{minipage}
\end{figure}

Let us end this section by considering the attractor of the map in Figure~\ref{fig:map}. The dimension of the attractor is
\[
  \dimH \Lambda = 1.488 \ldots
\]
There is a picture of the attractor $\Lambda$ in Figure~\ref{fig:att}. 

We may also consider the dimension of $\Lambda$ when $\gamma =1.8$, $\lambda = 0.5$ and $k = 0.1$. Then
\[
  \dimH \Lambda = 1.848\ldots
\]
A picture of this attractor is in Figure~\ref{fig:att2}. Both pictures where drawn by calculating the iterates of a small curve with tangents in the unstable cones.

\section{Proof of Theorem \ref{the:dimension}} \label{sec:proof}

Assume that $f$ satisfies condition \eqref{T} with $(\varepsilon_0, \delta)$-intersections. Let $\varepsilon > 0$.

\subsection{Coding of the system}

Let $\hat{f}$ be the lift of $f$ as described in Section~\ref{ssec:noninv} and let $\hat{\Lambda}$ denote the attractor of $\hat{f}$. We start by introducing a coding of the system $\hat{f} \colon \hat{\Lambda} \to \hat{\Lambda}$. If $\hat{x} \in \hat{\Lambda}$ then there is a sequence $\seq{\hat{s}} (\hat{x}) = \{ i_k\}_{k\in \Z}$ such that $\hat{f}^k (\hat{x}) \in \hat{K}_{i_k}$ for every $k\in \Z$. We let $\Sigma = \Sigma (\hat{\Lambda})$ be the set of all such sequences, that is $\Sigma (\hat{\Lambda}) = \seq{\hat{s}} (\hat{\Lambda})$.
Then there is an one-to-one correspondance $\rho \colon \Sigma \to \hat{D} \subset \hat{\Lambda}$, defined in the natural way. Let $\pi \colon \hat{K} \to K$ be the projection $\pi (x,y) = x$. 
Given a sequence $\seq{a} = \{ a_i \}_{i\in\Z}$, we define the cylinder set $ {}_{k} [ \seq{a} ]_l $ by
\[
  {}_{k} [ \seq{a} ]_l := \{\, \seq{b} = \{ b_i \}_{i\in\Z} \in \Sigma : b_i = a_i, \ \forall i=k,k+1,\ldots, l \,\}.
\]
The sets $\rho ( {}_{k} [ \seq{a} ]_l )$ and $\pi ( \rho ( {}_{k} [ \seq{a} ]_l ) )$ will also be called cylinders.

\subsection{Images of curves}

In this section we make use of condition \eqref{A9}, that the multiplicity entropy is zero, to get some estimates.

For $r \in \N$, we let $\mathcal{D}_r (\varepsilon)$ be the set of $r$-cylinders ${}_0 [\seq{a} ]_{r-1}$ such that there exists a point $p \in {}_0 [\seq{a} ]_{r-1}$ with
\begin{align*}
  e^{(\chi_\unstable - \varepsilon) r} & \leq \norm{ \diff{}_p (f^r) } \leq e^{(\chi_\unstable + \varepsilon) r}, & \forall v \in C^\unstable (p), \\
  e^{(\chi_\stable - \varepsilon) r} & \leq \norm{ \diff{}_p (f^r) } \leq e^{(\chi_\stable + \varepsilon) r}, & \forall v \in C^\stable (p).
\end{align*}

Let $q, r \in \N$, $l > 0$ and let $\gamma$ be a curve of length $l$ with tangents in the unstable cones. Let $\mathcal{W}_0 =\{\gamma\}$. We define $\mathcal{W}_n$ inductively. If $\mathcal{W}_{n-1}$ is a collection of curves, then we let $\mathcal{W}_n$ be the set of curves that are connected pieces of length between $l$ and $2l$, contained in the union of $\mathcal{D}_r$ and in some $f^q (\sigma)$, $\sigma \in \mathcal{W}_{n-1}$.

Since we require that the length of the curves in $\mathcal{W}_n$ are between $l$ and $2l$, the set $\mathcal{W}_n$ might not be uniquely defined, since there are several ways to divide a curve of length larger than $2l$ into pieces of length between $l$ and $2l$. It is however not important how this is done, so we will not give a precise definition of $\mathcal{W}_n$.

\begin{lemma} \label{lem:cutandexpand}
  Let $f\colon K \setminus N \to K$ satisfy the conditions \eqref{A2}, \eqref{A4}, \eqref{A5}--\eqref{A9}. For any $\varepsilon > 0$, there exist constants $C$, $q$, $r$ and $l > 0$, and a curve $\gamma$ with tangents in the unstable cones, such that if $N(n)$ denotes the number of curves in $\mathcal{W}_n$, then
  \[
    C^{-1} e^{(\chi_\unstable - \varepsilon) q (n - k)} \leq \frac{N(n)}{N(k)} \leq C e^{(\chi_\unstable + \varepsilon) q (n - k)},
  \]
  holds for all $n \geq k \geq 1$, and the derivatives of $f^q$ at a point $p \in W \in \mathcal{W}_n$ satisfies
  \begin{align}
    C^{-1} e^{(\chi_\unstable - \varepsilon) q k} \norm{v} &\leq \norm{ \diff{}_p (f^{qk}) (v) } \leq C e^{(\chi_\unstable + \varepsilon) q k} \norm{v} , & \forall v \in C^\unstable (p), \label{chiu} \\
    C^{-1} e^{(\chi_\stable - \varepsilon) q k} \norm{v} &\leq \norm{ \diff{}_p (f^{qk}) (v) } \leq C e^{(\chi_\stable + \varepsilon) q k} \norm{v} , & \forall v \in C^\stable (p). \label{chis}
  \end{align}

\end{lemma}

\begin{proof}
  Since the multiplicity entropy is zero, we can take $q$ large and $l > 0$ small, so that any curve of length $l$ with tangents in the unstable cone is cut in at most $e^{\varepsilon q}$ pieces when mapped by $f^q$.

  Since the Lebesgue measure of the complement of the union of $\mathcal{D}_r$ vanishes as $r \to \infty$, we can choose $r$ large so that the Lebesgue measure of the union of $\mathcal{D}_r$ is as close to that of $K$ as we like. Using property \eqref{A6}, we see that it is even possible to choose $r$ so large that the intersection of the complement of the union of $\mathcal{D}_r$ with any curve of length at least $l$ with tangents in the unstable cone has as small one dimensional Lebesgue measure as we like.

  Hence by first choosing $q$ and $l$, and then $r$ depending on $l$, it is possible to achieve that the sums of the lengths of the curves in $\mathcal{W}_n$ satisfies
  \[
    C_0 e^{(\chi_\unstable - \varepsilon) q n} \leq \sum_{\sigma \in \mathcal{W}_n} \mathrm{length} \, (\sigma) \leq C_0 e^{(\chi_\unstable + \varepsilon) q n},
  \]
  where $C_0$ is a constant depending on $f$, $q$, $l$ and $r$. This implies that the number of curves in $\mathcal{W}_n$ satisfies the statement in the lemma.
\end{proof}

\subsection{Frostman's lemma}

We define a probability measure $\mu_n$ with support on $\cup \mathcal{W}_n$ by
\[
  \mu_n = \frac{1}{N(n)} \sum_{W \in \mathcal{W}_n} \nu_W,
\]
where $\nu_W$ denotes the normalised Lebesgue measure on the curve $W$, and $N (n)$ denotes the number of elements in $\mathcal{W}_n$ as in Lemma~\ref{lem:cutandexpand}.

By taking a subsequence we can achieve that $\mu_{n}$ converges weakly to a probability measure $\mu$ with support in $\Lambda$. This measure will not be invariant, but its conditional measures on unstable manifold will be absolutely continuous with respect to the corresponding conditional measures of the \srb-measure, almost surely.

We will use the following method, originating from Frostman \cite{Frostman}, to estimate the dimension of $\Lambda$. If
\[
  \int \! \! \int \frac{\diff{\mu} (x) \diff{\mu} (y)}{|x-y|^s} < \infty,
\]
then $\dimH \Lambda \geq \dimH \supp \mu \geq s$. For a proof of this, see Falconer's book \cite{Falconer}.

Let $M$ be a number. Then
\begin{multline*}
  \int \! \! \int \min \biggl\{ M,\, \frac{1}{|x-y|^s} \biggr\} \, \diff{\mu_{n}} (x) \diff{\mu_{n}} (y) \\ \to \int \! \! \int \min \biggl\{ M,\, \frac{1}{|x-y|^s} \biggr\} \, \diff{\mu} (x) \diff{\mu} (y), \quad \mathrm{as }\ n \to \infty,
\end{multline*}
and
\begin{multline*}
  \int \! \! \int \min \biggl\{ M,\, \frac{1}{|x-y|^s} \biggr\} \, \diff{\mu} (x) \diff{\mu} (y)\\ \to \int \! \! \int \frac{1}{|x-y|^s} \, \diff{\mu} (x) \diff{\mu} (y), \quad \text{as } M \to \infty.
\end{multline*}

We will therefore estimate
\[
  E_s (n,M) = \int \! \! \int \min \biggl\{ M,\, \frac{1}{|x-y|^s} \biggr\} \, \diff{\mu_{n}} (x) \diff{\mu_{n}} (y).
\]
It is clear that $E_s (n,M) \leq M$. By the definition of the measure $\mu_n$ we immediately get that
\begin{equation} \label{sumint}
 E_s (n,M) = \sum_{W, V \in \mathcal{W}_{n}} \frac{1}{N (n)^2} \int \! \! \int \min \biggl\{ M,\, \frac{1}{|x-y|^s} \biggr\} \, \diff{\nu_{V}} (x) \diff{\nu_{W}} (y).
\end{equation}

We rewrite \eqref{sumint} as
\[
  E_s (n) = J_1 + J_2,
\]
with
\begin{align*}
J_1 &= \sum_{W \in \mathcal{W}_{n}} \frac{1}{N(n)^2} \int \! \! \int \min \biggl\{ M,\, \frac{1}{|x-y|^s} \biggr\} \, \diff{\nu_{W}} (x) \diff{\nu_{W}} (y),\\
J_2 &= \sum_{\substack{W, V \in \mathcal{W}_{n}, \\ V \ne W}} \frac{1}{N(n)^2} \int \! \! \int \min \biggl\{ M,\, \frac{1}{|x-y|^s} \biggr\} \, \diff{\nu_{V}} (x) \diff{\nu_{W}} (y).
\end{align*}

To estimate $J_1$ we note that
\[
 \int \! \! \int \min \biggl\{ M,\, \frac{1}{|x-y|^s} \biggr\} \, \diff{\nu_{W}} (x) \diff{\nu_{W}} (y) \leq {M}.
\]
Hence
\[
  J_1 \leq \sum_{W \in \mathcal{W}_{n}} \frac{M}{N(n)^2} = \frac{M}{N(n)},
\]
and so $J_1 \to 0$ as $n \to \infty$.

We will now estimate $J_2$ and show that $J_2$ is bounded as $n \to \infty$, provided that $s$ is sufficiently small.

Let $m < n$ and $W \in \mathcal{W}_n$. Then there is a unique $\alpha \in \mathcal{W}_{n-m}$ such that $W \subset f^{qm} (\alpha)$. Let $W_{-m}$ denote the set $W_{-m} \subset \alpha$ such that $W = f^{qm} (W_{-m})$.

\begin{figure}
\setlength{\unitlength}{0.007\textwidth}
 \begin{center}
  \begin{picture}(130,100)
    \put(55,84){\small $ \alpha $}
    \put(47,59){\small $ \beta $}
    \put(50,28){\small $V_{-m}$}
    \put(35,14){\small $W_{-m}$}
    \put(62,36){\small $B_{\varepsilon_0}$}
    \put(0,0){\includegraphics[height=100\unitlength]{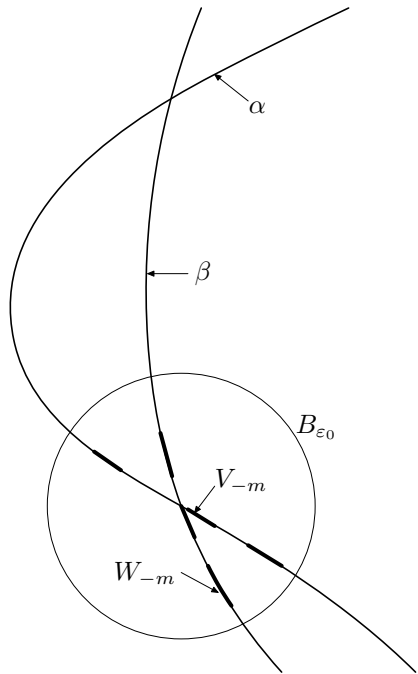}}
  \end{picture}
  
  \caption{A picture of intersections of unstable manifolds. The pre-images $W_{-m}$ and $V_{-m}$ are the thicker segments.} \label{fig:intersections}
 \end{center}
\end{figure}

Fix $m < n$ and take two different $\alpha$ and $\beta$ in $\mathcal{W}_{n-m}$ such that $\alpha_{-1}$ and $\beta_{-1}$ are in different cylinders. By condition \eqref{T} this implies that $\alpha$ and $\beta$ intersect $(\varepsilon_0, \delta)$-transversal. We will consider all manifolds $W$ and $V$ in $\mathcal{W}_n$ such that $W_{-m} \subset \alpha$, $V_{-m} \subset \beta$, and $W_{-m}$ and $V_{-m}$ are in the same $qm$-cylinder, which we denote by $S_m (W_{-m})$. There is a picture of this in Figure~\ref{fig:intersections}.

Note that $W$ and $V$ intersect if and only if $W_{-m}$ and $V_{-m}$ intersect, since $W_{-m}$ and $V_{-m}$ are in the same $qm$-cylinder.
If $W_{-m} \subset \alpha$ intersect $\beta$, then we estimate that
\begin{equation} \label{curveestimate}
  \sum_{\substack{V \in \mathcal{W}_n \\ V_{-m} \subset \beta \cap S_m (W_{-m}) }} \int \! \! \int \frac{1}{|x-y|^s} \, \diff{\nu_{V}} (x) \diff{\nu_{W}} (y) \leq C_1 e^{(\chi_\unstable - \chi_\stable + 2 \varepsilon)(s - 1) m},
\end{equation}
where $C_1$ does not depend on $W$, $\alpha$ and $\beta$.
Indeed, if $m$ is large, then we may assume that $W_{-m}$ and $V_{-m} \subset \beta \cap S_m (W_{-m})$ are contained in a ball of radius $\varepsilon_0$, and so the manifolds $f^{qm} (\beta)$ and $W$ intersect $(\varepsilon_0, C^2 e^{(\chi_\unstable - \chi_\stable + 2\varepsilon)m} \delta)$-transversal and we can estimate that
\begin{align*}
  \sum_{\substack{V \in \mathcal{W}_n \\ V_{-m} \subset \beta \cap S_m (W_{-m})}} \int \! \! \int \frac{1}{|x-y|^s} \, \diff{\nu_{V}} (x) \diff{\nu_{W}} (y) \leq C_0
  \int_{\gamma_1} \! \! \int_{\gamma_2} \frac{1}{|x-y|^s} \, \diff{} x \diff{} y,
\end{align*}
where $\gamma_1$ and $\gamma_2$ are the curves
\begin{align*}
  \gamma_1 &= \{\, (x_1, x_2) : x_1 = 0, |x_2| < \diam K \,\},\\
  \gamma_2 &= \{\, (x_1, x_2) : |x_2| < l, x_2 = C^2 e^{(\chi_\unstable - \chi_\stable + 2 \varepsilon)m} {\delta} x_1 \,\},
\end{align*}
and $C_0$ is a constant, that depends only on the second derivative of the map and the constants $\diam K$ and $l$.
To prove \eqref{curveestimate}, one easily checks that there exists a constant $C_1$ such that
\[
  C_0 \int_{\gamma_1} \! \! \int_{\gamma_2} \frac{1}{|x-y|^s} \, \diff{} x \diff{} y \leq C_1 e^{(\chi_\unstable - \chi_\stable + 2 \varepsilon)(s-1)m}.
\]

We consider now those manifolds $W$, such that $W_{-m} \subset \alpha$ does not intersect $\beta$. First, we consider those $V$ such that $W_{-m}$ and $V_{-m}$ lies in some ball $B_{\varepsilon_0}$ in the spirit of \eqref{T}.
If the distance between $W_{-m}$ and $\beta$ is $d (W_{-m}, \beta)$, then the distance between $W$ and $V \subset f^{qm} (\beta)$ is larger than $C^{-1} e^{(\chi_s - \varepsilon) m} d (W_{-m}, \beta) $ by \eqref{chis}. If we choose the length $l$ in the construction of $\mathcal{W}_n$ sufficiently small, then 
we can approximate the integral by 
\[
  \sum_{\substack{V \in \mathcal{W}_n \\ V_{-m} \subset \beta \cap S_m (W_{-m}) }} \int \! \! \int \frac{1}{|x-y|^s} \, \diff{\nu_{V}} (x) \diff{\nu_{W}} (y) \leq l^{-2} \int_{\gamma_1} \! \! \int_{\gamma_2} \frac{1}{|x-y|^s} \, \diff{} x \diff{} y,
\]
where $\gamma_1$ and $\gamma_2$ are two parallell line segments of length $l$ and $\diam K$, and with distance $d (W_{-m}, \beta) / 2$. The last integral is estimated by
\begin{multline*}
  \int_{\gamma_1} \! \! \int_{\gamma_2} \frac{1}{|x-y|^s} \, \diff{} x \diff{} y 
  \leq \int_{-\infty}^\infty \frac{1}{\bigl( \sqrt{x^2 + (d (W_{-m}, \beta) / 2)^2 } \bigr)^s} \, \diff{} x \\
   = 2^{s} e^{(\chi_s - \varepsilon) (1 - s) m} d (W_{-m}, \beta)^{1-s} \int_0^\infty \frac{\mathrm{d}x}{(1 + x^2)^{\frac{s}{2}}},
\end{multline*}
and so
\begin{multline}
  \sum_{\substack{V \in \mathcal{W}_n \\ V_{-m} \subset \beta \cap S_m (W_{-m}) }} 
  \int \! \! \int \frac{1}{|x-y|^s} \, \diff{\nu_{V}} (x) \diff{\nu_{W}} (y) \\
  \leq C_2 e^{(\chi_s - \varepsilon) (1 - s) m} d (W_{-m}, \beta)^{1-s}, \label{dest}
\end{multline}
for some constant $C_2$, provided that $s>1$. In fact, one easily shows that
\[
  \int_0^\infty \frac{\mathrm{d}x}{(1+x^2)^{\frac{s}{2}}}
  = \frac{\sqrt{\pi}}{2} \frac{ \Gamma(\frac{s-1}{2} )}{ \Gamma (\frac{s}{2})},
\]
by the change of variable $t = \frac{1}{1+x^2}$ and the observation that
\[
  \frac{\sqrt{\pi}}{2} \frac{\Gamma(\frac{s-1}{2} )}{ \Gamma (\frac{s}{2})}
  = \frac{1}{2} \frac{\Gamma (\frac{1}{2}) \Gamma(\frac{s-1}{2})}{ \Gamma (\frac{s}{2})}
  = \frac{1}{2} B( {\textstyle \frac{1}{2}, \frac{s-1}{2} } ),
\]
where $B(x,y) = \frac{ \Gamma(x) \Gamma(y)}{\Gamma ( x +y )} = \int_0^\infty t^{x-1} (1-t)^{y-1}\, \mathrm{d} t$ is the beta function.

We cover the intersections of $\alpha$ and $\beta$ by balls $B_{\varepsilon_0}$. Since $K$ is a bounded set the number of such balls will always be less than some number $N_B$. Given a manifold $W_{-m} \subset \alpha$, either $W_{-m}$ intersect one of these balls or lies a distance of at least $\varepsilon_0$ from each of the intersections of $\alpha$ and $\beta$.

If $W_{-m}$ lies in $B_{\varepsilon_0}$, with distance $d_{W_{-m}}$ from the center of the ball, then by property \eqref{T}, the distance between $W_{-m}$ and $V_{-m} \subset \beta$ is at least $\delta d_{W_{-m}}$. The manifolds $W_{-m}$ and $V_{-m}$ are subsets of the two larger manifolds $\alpha$ and $\beta$ (see Figure~\ref{fig:intersections}). On each side of the intersection of these larger manifolds (or the closest point in case they do not intersect) we can enumerate the pairs $W_{-m}$ and $V_{-m}$, such that the distance from the center of the ball $B_{\varepsilon_0}$ to the $i$th manifold $W_{-m}$ is increasing. Since two different $W_{-m}$ do not intersect, the distance from the center of $B_{\varepsilon_0}$ to the $i$th manifold $W_{-m}$ is at least $i \frac{l}{C e^{(\chi_\unstable + \varepsilon) m}}$, since the length of each $W_{-m}$ is at least $\frac{l}{C e^{(\chi_\unstable + \varepsilon) m}}$ by \eqref{chiu}. (We measure the distance along the large manifold containing all the $W_{-m}$.) 
This implies that the distance between the $i$th $W_{-m}$ and $V_{-m} \subset \beta$ is at least
$\delta i \frac{l}{C e^{(\chi_\unstable + \varepsilon) m}}$ and so the distance between the corresponding $W$ and $V$ is at least $C^{-2} e^{ (\chi_s - \varepsilon) m } \delta i \frac{l}{e^{(\chi_\unstable + \varepsilon) m }}$. There are at most $C_3 M (m)$ different $W_{-m}$ in $B_{\varepsilon_0}$, where $C_3$ is a constant that dependds on $l$ and $M (m)$ denotes the number of $qm$-cylinders.
By \eqref{curveestimate} and \eqref{dest} we estimate that
\begin{align}
  &\sum_{ \substack{ W_{-m} \subset \alpha \cap B_{\varepsilon_0} ,\\  V_{-m} \subset \beta \cap S_m (W_{-m}) \cap B_{\varepsilon_0} } }
  \int \! \! \int \frac{1}{|x - y|^s} \, \diff{\nu_{V}} (x) \diff{\nu_{W}} (y) \nonumber \\ 
  & \hspace{1cm} < C_1 e^{(\chi_\unstable - \chi_\stable + 2 \varepsilon)(s-1)m} + 2 \sum_{i=1}^{ C_3 M (m) } C_2 \Bigl( C^{-2} e^{ (\chi_\stable - \varepsilon) m } \delta i \frac{l_0}{e^{(\chi_\unstable + \varepsilon) m }} \Bigr)^{1-s} \nonumber \\
  & \hspace{1cm} < C_4 e^{(\chi_\unstable - \chi_\stable + 2\varepsilon)(s-1)m} M (m)^{2-s}. \label{inball}
\end{align}
If we sum over the balls $B_{\varepsilon_0}$ needed to cover the intersection of $\alpha$ and $\beta$, we get
\begin{multline} \label{est1}
   \sum_{B_{\varepsilon_0}}
   \sum_{ \substack{ W_{-m} \subset \alpha \cap B_{\varepsilon_0} ,\\  V_{-m} \subset \beta \cap S_m (W_{-m}) \cap B_{\varepsilon_0} } }
   \int \! \! \int \frac{1}{|x-y|^s} \, \diff{\nu_{V}} (x) \diff{\nu_{W}} (y) \\
   < C_5 e^{(\chi_\unstable - \chi_\stable + 2 \varepsilon)(s-1)m} M (m)^{2-s}.
\end{multline}

For those $W$ and $V$ such that $W_{-m}$ and $V_{-m}$ are not inside a ball $B_{\varepsilon_0}$ we have $d(W, V) > C^{-1} e^{( \chi_\stable - \varepsilon) m} \varepsilon_0$, and estimate by \eqref{dest} that
\begin{multline} \label{est2}
  \int \! \! \int \frac{1}{|x-y|^s} \, \diff{\nu_{V}} (x) \diff{\nu_{W}} (y) \\ \leq C_2 e^{( \chi_\stable - \varepsilon) (1-s) m} C^{s-1} \varepsilon_0^{1-s} = C_6 e^{( \chi_\stable - \varepsilon) (1-s) m}.
\end{multline}
The number of such pairs $W$ and $V$ are at most some constant $C_7$ times the number of $qm$-cylinders, denoted by $M (m)$. We get by \eqref{est1} and \eqref{est2} that
\begin{multline}
  \sum_{\substack{W,V \in \mathcal{W}_n,\\ W_{-m} \subset \alpha, \\ V_{-m} \subset \beta } }
  \int \! \! \int \frac{1}{|x-y|^s} \, \diff{\nu_{V}} (x) \diff{\nu_{W}} (y) \\
  < C_5 e^{(\chi_\unstable - \chi_\stable + 2 \varepsilon)(s-1)m}  M (m)^{2-s} + C_6 C_7 e^{(\chi_\unstable -\chi_\stable (s-1) + s \varepsilon) m} M (m). \label{outsideball}
\end{multline}

We will now sum over all $m$, $\alpha$ and $\beta$, and write $J_2$ as
$J_2 = J_3 + J_4$, with
\begin{align*}
  J_3 &= \sum_{m=0}^{n-1} \sum_{\substack{\alpha, \beta \in \mathcal{W}_{n-m} \\ \alpha \ne \beta }} \sum_{\substack{W,V \in \mathcal{W}_n,\\ W_{-m} \subset \alpha, \\ V_{-m} \subset \beta } } \frac{\int \! \! \int \min \biggl\{ M,\, \frac{1}{|x-y|^s} \biggr\} \, \diff{\nu_{V}} (x) \diff{\nu_{W}} (y)}{N(n)^2} , \\
  J_4 &= \sum_{m=0}^{n-1} \sum_{\alpha \in \mathcal{W}_{n-m} } \sum_{\substack{W,V \in \mathcal{W}_n,\\ W_{-m}, V_{-m} \subset \alpha, \\ W_{-m} \ne V_{-m} } } \frac{\int \! \! \int \min \biggl\{ M,\, \frac{1}{|x-y|^s} \biggr\} \, \diff{\nu_{V}} (x) \diff{\nu_{W}} (y)}{N(n)^2} .
\end{align*}
Similarly as for $J_1$ we obtain that $J_4 \to 0$ as $n \to \infty$. It remains to estimate $J_3$.

Using that there are $N(n-m)$ different $\alpha$ and $\beta$, we get by \eqref{inball} and \eqref{outsideball} that
\begin{multline*}
   J_3 \leq \sum_{m=0}^{n-1}  N (n - m)^2 \frac{ C_5 e^{(\chi_\unstable - \chi_\stable + 2 \varepsilon)(s - 1) q m}  M (m)^{2 - s}}{N(n)^2} \\ + \sum_{m = 0}^{n - 1}  N(n - m)^2 \frac{ C_6 C_7 e^{(\chi_\unstable - \chi_\stable (s - 1) + s \varepsilon) q m} M (m) }{N(n)^2}.
\end{multline*}
We now use that the topological entropy is $\chi_\unstable$, and thus the number of cylinders satisfy $M (m) \leq C_6 e^{(\chi_\unstable + \varepsilon) q m}$, for some constant $C_6$. This yields
\[
  J_3 \leq C_8 \sum_{m=0}^{n-1}  \frac{N (n-m)^2 e^{(\chi_\unstable - \chi_\stable (s-1) + s \varepsilon) q m} }{N(n)^2},
\]
for some constant $C_8$ that does not depend on $n$. By Lemma~\ref{lem:cutandexpand} we have $N (n-m) / N(n) \leq C e^{ - (\chi_\unstable - \varepsilon) q m}$, so
\[
  J_3 \leq C_8 \sum_{m=0}^{n-1} e^{ (-\chi_\unstable - \chi_\stable (s - 1) + s \varepsilon ) q m}.
\]

We conclude that $J_3$ is bounded as a function of $n$ provided that
$- \chi_\unstable - (s - 1) \chi_\stable + s \varepsilon < 0$ and $s < 2$ or equivalently
\begin{equation} \label{scondition}
  s < 1 - \frac{\chi_\unstable - \varepsilon}{ \chi_\stable - \varepsilon} \quad \mathrm{and} \quad s < 2.
\end{equation}

We have therefore obtained that, if $s$ satisfies \eqref{scondition}, then $J_1$ and $J_2 = J_3 + J_4$ are bounded, and so the integral
\[
  \int \! \! \int \min \biggl\{ M,\, \frac{1}{|x-y|^s} \biggr\} \, \diff{\mu} (x) \diff{\mu} (y)
\]
is uniformly bounded and hence converges as $M \to \infty$. This proves that
\[ \int \! \! \int \frac{1}{|x-y|^s} \, \diff{\mu} (x) \diff{\mu} (y) < \infty,
\]
provided that (\ref{scondition}) holds.
Hence
\[
  \dimH \Lambda \geq \min \biggl\{2,\, 1 - \frac{\chi_\unstable - \varepsilon}{ \chi_\stable - \varepsilon} \biggr\}.
\]
Let $\varepsilon \to 0$.

\end{document}